\documentclass[12pt]{article}
\usepackage{amsmath, amsthm, amssymb}
\newtheorem{theorem}{Theorem}[section]
\newtheorem{lemma}[theorem]{Lemma}
\theoremstyle{definition}
\newtheorem{definition}[theorem]{Definition}

\theoremstyle{remark}

\numberwithin{equation}{section}

\begin{document}
\title{Actions of $SL(n,\mathbb{Z})$ on homology spheres}
\author{Kamlesh Parwani}
%\address{Department of Mathematics, Northwestern University, Evanston IL 60208, USA}
%\email{forty2@math.northwestern.edu}
%\subjclass{Primary 57S25; Secondary 37C85, 57S17}
%\keywords{Homology spheres, group actions, almost simple}
\date{January 25, 2004.}

\maketitle

\begin{abstract}
Any continuous action of $SL(n,\mathbb{Z})$, where $n>2$, on a $r$%
-dimensional mod $2$ homology sphere factors through a finite group action
if $r<n-1$. In particular, any continuous action of $SL(n+2,\mathbb{Z})$ on
the $n$-dimensional sphere factors through a finite group action.
\end{abstract}

%%%%%%%%%%%%%%%%%%%%%%%%%%%%%%%%%%%%%%%%%%%%%%%%%%%%%%%%%%%%%%%%%%%%%

\section{Introduction}

In this paper, we study the actions of $SL(n,\mathbb{Z})$ on spheres and,
more generally, actions on homology spheres. The group $SL(n,\mathbb{Z})$
acts on the $(n-1)$-dimensional sphere via the linear action on vectors in $%
\mathbb{R}^{n}$. This action is minimal in the following sense.

\begin{theorem}
Any continuous action of $SL(n,\mathbb{Z})$, where $n>2$, on a $r$%
-dimensional mod $2$ homology sphere factors through a finite group action
if $r<n-1$.
\end{theorem}

Since all spheres are mod 2 homology spheres, any continuous action of $%
SL(n+2,\mathbb{Z})$ on a $n$-dimensional sphere factors through a finite
group action. This result supports the following conjecture of Farb and
Shalen (see \cite{Farb}). \smallskip

\hspace{-14 pt} \textbf{Conjecture}. \textit{Any smooth action of a
finite-index subgroup of $SL(n,\mathbb{Z})$, where $n>2$, on a $r$%
-dimensional compact manifold factors through a finite group action if $%
r<n-1 $.} \smallskip

This conjecture is an analogue of a special case of one of the central
conjectures in the Zimmer program (see \cite{Zimmer}). Theorem 1.1 may also
be viewed as a (partial) generalization of Witte's theorem in \cite{Witte}.

\begin{theorem}[Witte]
If $\Gamma $ is a subgroup of finite index in $SL(n,\mathbb{Z})$ with $n\geq
3$, then every continuous action of $\Gamma $ on the circle factors through
a finite group action.
\end{theorem}

We cannot obtain our result for finite-index subgroups because we rely
heavily on the existence of finite order elements in $SL(n,\mathbb{Z})$, and
there are subgroups of finite index in $SL(n,\mathbb{Z})$ that have no
elements of finite order (see Corollary 6.13 in \cite{Rg}).

The paper is organized in the following manner. In section 2, we prove the
existence of certain desirable finite-order elements in $SL(n,\mathbb{Z})$
and reduce the problem to a problem of a finite group action. In section 3,
we use some classical results from the theory of compact transformation
groups to prove that certain groups cannot act effectively (faithfully) on
homology spheres and show that these results imply Theorem 1.1.

In section 4 we observe that the action of $SL(n,\mathbb{Z})$ is trivial on
low-dimensional spheres. This result is analogous to the following theorem
by Weinberger (see \cite{Shmuel}).

\begin{theorem}[Weinberger]
The discrete group $SL(n,\mathbb{Z})$, with $n\geq 3$, can act smoothly on
the torus $T^{m}$, $m<n$, only trivially.
\end{theorem}

\subsection*{Definitions and Notation}

For us a $n$-dimensional mod $2$ homology sphere is a locally compact,
finite-dimensional Hausdorff space $X$ such that $H^{*}(X;\mathbb{Z}%
_{2})=H^{*}(S^{n};\mathbb{Z}_{2})$, where $H^{*}$ means cohomology with
compact support and $S^{n}$ is the $n$-dimensional sphere. Furthermore, the $%
n$-dimensional mod 2 homology sphere is also a $n$-dimensional generalized
manifold---a finite-dimensional metric \textbf{ANR} with $H^{n}(X,X-x;%
\mathbb{Z})=\mathbb{Z}$ and $H^{i}(X,X-x;\mathbb{Z})=0$ for all $i\neq n$,
for all points $x\in X$ (see \cite{Bredon2} for equivalent conditions).

In section 2, we will prove the existence of certain finite groups in $SL(n,%
\mathbb{Z})$; for this we need manageable notation to deal with matrices.
Consider the standard basis of vectors in $\mathbb{R}^{n}$ and let $%
(1)=(1,0,0,..,0)$, $(2)=(0,1,0,...,0)$,..., $(n)=(0,0,0,...,1)$, $%
(-1)=(-1,0,0,..,0)$, $(-2)=(0,-1,0,...,0)$,..., and $(-n)=(0,0,0,...,-1)$.
The permutation $(i,j)$, with $1\leq i,j\leq n$, is the matrix that takes
the vector $(i)$ to the vector $(j)$ and the vector $(j)$ to $(i)$. For
example, we have the following matrices in $SL(4,\mathbb{Z})$. 
\begin{equation*}
(1,-1)(2,-2)=\left( 
\begin{array}{cccc}
-1 & 0 & 0 & 0 \\ 
0 & -1 & 0 & 0 \\ 
0 & 0 & 1 & 0 \\ 
0 & 0 & 0 & 1
\end{array}
\right) ,(1,2)(3,4)=\left( 
\begin{array}{cccc}
0 & 1 & 0 & 0 \\ 
1 & 0 & 0 & 0 \\ 
0 & 0 & 0 & 1 \\ 
0 & 0 & 1 & 0
\end{array}
\right)
\end{equation*}

This notation is unconventional, but it will simplify matrix multiplication
by reducing it to the standard procedure of multiplying permutations. Also,
let $I$ and $-I$ be the identity matrix and the negative identity matrix
respectively.

%%%%%%%%%%%%%%%%%%%%%%%%%%%%%%%%%%%%%%%%%%%%%%%%%%%%%%%%%%%%%%%%%%%%%%%%

\section{Almost simple groups and $SL(n,\mathbb{Z})$}

\begin{definition}
An element $g$ in a group $G$ is \textit{central} if $g$ commutes with every element
in $G$. A subgroup $H$ of the group $G$ is central if every element of $H$
is central. Also let $Z$ denote the \textit{center}, the subgroup of all central elements.
\end{definition}

It will be clear from the context what is meant by $Z$. For example, $G_1/Z$ and $G_2/Z$ are factor groups obtained when the groups $G_1$ and $G_2$ are quotiented by their respective centers. 

\begin{definition}
A group $G$ is \textit{almost simple} if every normal subgroup is either
finite and central, or has finite index in $G$.
\end{definition}

The Margulis normal subgroups theorem (see \cite{Margulis}) asserts that an
irreducible lattice in a semi-simple Lie group with $\mathbb{R}$-rank $\geq
2 $ is almost simple. In particular, $SL(n,\mathbb{Z})$ is almost simple for $n\geq
3$.

The following lemma follows easily from Margulis' Theorem and the definition
of an almost simple group.

\begin{lemma}
Let $\phi :SL(n,\mathbb{Z})\rightarrow H$ be a homomorphism where $n\geq 3$.
If $\phi (g)=1$ for some non-central element $g$, then the kernel of $\phi $
is a finite-index subgroup, and therefore, $\phi $ factors through a
homomorphism of a finite group.
\end{lemma}

So to prove Theorem 1.1, it suffices to show that a
finite-order, non-central element acts trivially. In the next section we
prove that there always exists an involution which acts trivially. Now we
show that subgroups containing non-central involutions always exist.

\begin{theorem}
If $n$ is odd, there exists a subgroup of $SL(n,\mathbb{Z})/Z$
isomorphic to $(\mathbb{Z}_{2})^{n-1}$. If $n$ is even, there
exists a subgroup of $SL(n,\mathbb{Z})/Z$ isomorphic to $(\mathbb{Z}_{2})^{n-2}$.
\end{theorem}

\begin{proof}
First note that $Z$ is just the identity matrix when $n$ is odd and only
contains the identity and the negative identity matrices when $n$ is even.
The theorem is very easy to prove when $n$ is odd. Simply count the elements
in the subgroup of diagonal matrices that have only $\pm 1$ in the diagonal
entries. This number is given by the following formula (which is true for
even and odd $n$).

\begin{equation*}
\left( 
\begin{array}{c}
n \\ 
0
\end{array}
\right) +\left( 
\begin{array}{c}
n \\ 
2
\end{array}
\right) +\left( 
\begin{array}{c}
n \\ 
4
\end{array}
\right) +....=2^{n-1}
\end{equation*}

So we have an abelian subgroup containing $2^{n-1}$ elements, all of order $2$. By the classification of finite abelian groups, this group is
isomorphic to $(\mathbb{Z}_{2})^{n-1}$. This proves that there is always
exists a subgroup of $SL(n,\mathbb{Z})$ isomorphic to $(\mathbb{Z}_{2})^{n-1} $.

When $n$ is even, there are $2^{n-2}$ cosets of diagonal matrices. So we
have a subgroup isomorphic to $(\mathbb{Z}_{2})^{n-2}$ in $SL(n,\mathbb{Z})/Z $.
\end{proof}

\begin{lemma}
There exists a subgroup of $SL(4,\mathbb{Z})/Z$ isomorphic to $(\mathbb{Z}%
_{2})^{3}$.
\end{lemma}
\begin{proof}
The above theorem states that there are 4 cosets of diagonal matrices. We
observe that the equivalence class of $(1,2)(3,4)$ commutes with the
equivalence classes of diagonal matrices in $SL(4,\mathbb{Z})/Z$.

$(1,-1)(2,-2)(1,2)(3,4)=(1,-2)(3,4)=(1,2)(3,4)(1,-1)(2,-2)$

\begin{multline*}
(1,-1)(3,-3)(1,2)(3,4)=(1,2,-1,-2)(3,4,-3,-4)= \\
-I(1,-2,-1,2)(3,-4,-3,4)=(1,2)(3,4)(1,-1)(3,-3)
\end{multline*}
All other cases are similar or easier. So we have a subgroup isomorphic to $(\mathbb{Z}_{2})^{3}$ in $SL(4,\mathbb{Z})/Z$.
\end{proof}

Note that for any action of $SL(n,\mathbb{Z})$ or $SL(n,\mathbb{Z})/Z$ on an
orientable, generalized manifold $M$, the induced action of the subgroups in Theorem 2.4 and Lemma 2.5 on $M$ is orientation-preserving. This is because $SL(n,\mathbb{Z})$ is a perfect group, that is, it is equal to its commutator subgroup. Any homomorphism from a perfect group to an abelian group must be trivial. In particular, any homomorphism from $SL(n,\mathbb{Z})$ to $\mathbb{Z}_2$ must be trivial, and so, any action of $SL(n,\mathbb{Z})$ on an orientable, generalized manifold must be orientation-preserving. Then the same result holds for subgroups of $SL(n,\mathbb{Z})$ and $SL(n,\mathbb{Z})/Z$.

We now introduce new notation which will greatly simplify the proof and the
statement of the next lemma. In $SL(n,\mathbb{Z})$, let $%
-I_{i}^{j}=(i,-i)(i+1,-i-1)...(j,-j)$, where $1\leq i<j\leq n$ and $j-i$ is
odd. So when $n$ is even, $-I=-I_{1}^{n}$ and the center $Z=\left\langle
-I_{1}^{n}\right\rangle $.

\begin{lemma}
Consider the action of $SL(n,\mathbb{Z})/\left\langle
-I_{1}^{n}\right\rangle $, for even $n\geq 4$, on an orientable generalized
manifold $M$ and suppose that $-I_{n-1}^{n}$ acts on $M$ with a non-empty
fixed point set. Then there is an induced action of a subgroup isomorphic to 
$SL(n-2,\mathbb{Z})/\left\langle -I_{1}^{n-2}\right\rangle$ on the fixed
point set of $-I_{n-1}^{n}$.
\end{lemma}
\begin{proof}
Let $F$ be the fixed point set of $-I_{n-1}^{n}$. Now all elements of $SL(n,%
\mathbb{Z})$ that leave the vectors $(n-1)$ and $(n)$ fixed commute with $%
-I_{n-1}^{n}$. This implies that there is a subgroup isomorphic to $SL(n-2,%
\mathbb{Z})$ in the centralizer of $-I_{n-1}^{n}$, and therefore, it acts on 
$F$ modulo the action of $-I_{1}^{n}$. The action of $-I_{n-1}^{n}$ of $F$
is trivial, and so, the action of $-I_{1}^{n}$ on $F$ is equivalent to the
action of $-I_{1}^{n-2}$ on $F$. This proves that there is an action of $%
SL(n-2,\mathbb{Z})$ on $F$ modulo the action of $-I_{1}^{n-2}$; in other
words, there is an induced action of a subgroup isomorphic to $SL(n-2,%
\mathbb{Z})/\left\langle -I_{1}^{n-2}\right\rangle $ on the fixed point set
of $-I_{n-1}^{n}$.
\end{proof}

%%%%%%%%%%%%%%%%%%%%%%%%%%%%%%%%%%%%%%%%%%%%%%%%%%%%%%%%%%%%%%%%%%%%%%%%%%%%
\section{Actions on homology spheres}

We begin this section by recalling some very famous theorems from the theory
of compact transformation groups. Most of the results stated below are true
for all primes but we only consider the case when the prime is two. Good
references for these theorems are \cite{Bredon}, \cite{Borel}, \cite{Smith1}%
, and \cite{Smith2}.

\begin{theorem}[Smith]
The group $\mathbb{Z}_{2}\times \mathbb{Z}_{2}$ cannot act semifreely and
effectively on any mod $2$ homology sphere, that is, no two elements in the
group $\mathbb{Z}_{2}\times \mathbb{Z}_{2}$ have the same fixed-point set
(which may be empty) when the action is effective.
\end{theorem}

As a corollary, we get the result that the group $\mathbb{Z}_{2}\times %
\mathbb{Z}_{2}$ cannot act freely on any mod $2$ homology sphere.

\begin{theorem}[Smith]
When $\mathbb{Z}_{2}$ acts effectively on a $n$-dimensional mod $2$ homology
sphere $X$, the fixed set (pointwise) is a $m$-dimensional mod $2$ homology
sphere with $m<n$.
\end{theorem}

Furthermore since $X$ is a generalized manifold, one can distinguish between
actions that preserve orientation and actions that reverse orientation. If
orientation is preserved, $n-m$ is even (see \cite{Bredon2}).

The following theorem is due to Smith and is well known. In \cite{Smith2} Smith
proves that any action of $(\mathbb{Z}_{2})^{n+2}$ on a $n$-dimensional mod $%
2$ homology sphere cannot be effective. Theorem 3.3 follows from his proof
and other arguments in \cite{Smith1}. We provide a simple proof below.

\begin{theorem}
Any orientation-preserving action of $(\mathbb{Z}_{2})^{n+1}$ on a $n$%
-dimensional mod $2$ homology sphere cannot be effective.
\end{theorem}

\begin{proof}
We prove this theorem by induction. The result easily follows for the $0$%
-dimensional sphere---$S^{0}$. Every orientation-preserving action of $%
\mathbb{Z}_{2}$ on $S^{0}$ is trivial. Now assume that the theorem is true
for all $k\leq n-1$, that is, assume that any orientation-preserving action of $(\mathbb{Z}%
_{2})^{k+1}$ on a $k$-dimensional homology sphere is not effective for all $%
k\leq n-1.$

Suppose that an effective, orientation-preserving action of $(\mathbb{Z}_{2})^{n+1}$ on a $n
$ dimensional mod $2$ homology sphere exists. Theorem 3.1 asserts that this action cannot be free and so there
is some element acting with a non-empty fixed-point set. By theorem 3.2, we
know that this fixed set is always a $m$-dimensional mod $2$ homology sphere
with $m<n$. Let $h$ be the element with the property that the fixed-point
set of $h$, $F(h)$, is maximal among all other fixed sets of elements in $(%
\mathbb{Z}_{2})^{n+1}$, that is, $F(h)$ is not properly contained in the
fixed set of any other element in $(\mathbb{Z}_{2})^{n+1}$. Let $m$ be the
dimension of $F(h)$; note that $m\leq n-2$  since the action is orientation-preserving.

Let $G$ be the subgroup isomorphic to $(\mathbb{Z}_{2})^{n}$ such that $%
h\notin G$. Now because the action is abelian, there is an induced action of 
$G$ on $F(h)$, the $m$-dimensional mod $2$ homology sphere. This action of $G$ on
$F(h)$ may not be orientation-preserving, but in any case, there is an index 2 subgroup
isomorphic to $(\mathbb{Z}_{2})^{n-1}$ with an orientation-preseving action on $F(h)$.
The induction hypothesis implies that this action cannot be effective and so there exists
a $g\in G$ such that the fixed set of $g$ contains $F(h)$. The maximality of 
$F(h)$ implies that $h$ and $g$ have the same fixed-point set and this
contradicts Theorem 3.1. So any orientation-preserving action of $(\mathbb{Z}%
_{2})^{n+1}$ on a $n$-dimensional mod $2$ homology sphere cannot be
effective.
\end{proof}

Note that Theorem 3.3 is sharp, that is, there are effective actions of $(%
\mathbb{Z}_{2})^{n}$ on $n$-dimensional spheres that preserve orientation.
Theorem 2.4 shows that there exists a subgroup of $SL(n+1,\mathbb{Z})$
isomorphic to $(\mathbb{Z}_{2})^{n}$ and $SL(n+1,\mathbb{Z})$ acts
effectively on the $n$-dimensional sphere via the linear action on vectors
in $\mathbb{R}^{n+1}$. This effective action of $SL(n+1,\mathbb{Z})$ is
orientation-preserving and so the action of its subgroup isomorphic to $(%
\mathbb{Z}_{2})^{n}$ is also effective and preserves orientation.

We are now ready to prove the main theorem.

\begin{proof}[Proof of Theorem 1.1]
The theorem is very easy to prove in the case
when $r=0$ because the $0$-dimensional mod $2$ homology sphere is just the
set of two points---$S^{0}$. Now because the group $SL(n,\mathbb{Z})$ is perfect, it has an orientation-preserving action on $S^{0}$, and all orientation-preserving actions on $S^{0}$ are trivial.

So assume that $r\geq 1$. There are two cases to consider, $-I$ does not act
trivially and when $-I$ acts trivially.
\smallskip

CASE 1: $-I$ does not act trivially.

In this case, Theorem 2.4 asserts that there is a subgroup of $SL(n,\mathbb{Z})$ isomorphic to $(\mathbb{Z}_{2})^{n-1}$, and then it follows from general principles that the action of this subgroup is orientation-preserving. Theorem 3.3 implies that $(\mathbb{Z}_{2})^{n-1}$ cannot act effectively on any $r$-dimensional mod $2$ homology sphere, where $r<n-1$. So a non-central involution acts trivially and by the remark made after Lemma 2.3, this is sufficient to prove that the action of $SL(n,\mathbb{Z})$ on any $r$-dimensional mod $2$ homology sphere,
where $r<n-1$, factors through a finite group action. Note that if $n$ is
odd, $-I$ does not belong to $SL(n,\mathbb{Z})$, and so, the theorem has
been established for odd integers greater than 2.
\smallskip

CASE 2: $-I$ acts trivially.

When $-I$ acts trivially, $n$ must be even and the action factors through an
action of $SL(n,\mathbb{Z})/Z$. In this case, rather than showing only that some
non-central element of $SL(n,\mathbb{Z})$ acts trivially, we will prove the more precise
conclusion that some diagonal element conjugate to $-I_{1}^{2}$ acts trivially.

For $n=4$, Lemma 2.5 states that there is a subgroup isomorphic to $(\mathbb{Z}_{2})^{3}$ in $SL(4,\mathbb{Z})/Z$ and the result follows in a similar fashion. So every action of $SL(4,\mathbb{Z})$ on a $r$-dimensional mod 2 homology sphere, with $r<3$, factors through an action of a finite group. In the next section, we will prove that the action of $SL(4,\mathbb{Z})$ on low-dimensional mod 2 homology spheres---dimension less than 3---is trivial.
This implies that the action of the subgroup isomorphic to $(\mathbb{Z}_{2})^{3}$ in $SL(4,\mathbb{Z})/Z$ on mod 2 homology spheres with dimension less than 3 must be trivial. So the equivalence classes of diagonal elements conjugate to $-I_{1}^{2}$ act trivially; call these elements diagonal elements of length 2. In particular, for every action of $SL(4,\mathbb{Z})/Z$ on a mod 2 homology spheres with dimension less than 3, there exists an
equivalence class of a diagonal element of length 2 which acts trivially.

Now suppose that for every action of $SL(2k-2,\mathbb{Z})/Z$ on a mod 2
homology spheres with dimension less than $2k-3$, there exists an
equivalence class of a diagonal element of length 2 which acts trivially. Consider the action of $SL(2k,%
\mathbb{Z})/Z$, for $k>2$, on a $r$-dimensional mod $2$ homology sphere, where $r<2k-1$, and suppose that no class of diagonal elements
of length 2 acts trivially. The action of the equivalence classes of a
diagonal elements of length 2 cannot be free, by the corollary to Theorem
3.1, and so, without loss of generality assume that the fixed point set of
the class represented by $-I_{n-1}^{n}$ is maximal among all the fixed point
sets of equivalence classes of a diagonal elements of length 2. This fixed
point set must be a mod 2 homology sphere with dimension $m$, where $m<2k-3$,
since the action of $-I_{n-1}^{n}$ is orientation-preserving. Lemma 2.6
implies that we have an action of $SL(2k-2,\mathbb{Z})/Z$ on this $m$-dimensional
mod $2$ homology sphere and it now follows that there is an
equivalence class of a diagonal element of length 2 which acts trivially on
this $m$-dimensional mod $2$ homology sphere. The maximality of the fixed point set
of the class of $-I_{n-1}^{n}$ implies that two elements have identical fixed point sets and this contradicts Theorem 3.1. Therefore, there exists an equivalence class of a diagonal element of length 2 in $SL(n,\mathbb{Z})/Z$ that acts trivially, and it follows that some diagonal element conjugate to $-I_{1}^{2}$ in $SL(n,\mathbb{Z})$ acts trivially.

So when $n$ is even and greater than 3, there exists a non-central involution which acts trivially. By the remark made after Lemma 2.3, this is sufficient to prove that the action of $SL(n,\mathbb{Z})$ on a $r$-dimensional mod $2$ homology sphere, where $r<n-1$,
factors through an action of a finite group.
\end{proof} 
%%%%%%%%%%%%%%%%%%%%%%%%%%%%%%%%%%%%%%%%%%%%%%%%%%%%%%%%%%%%%%%%%%%%%%%%%%

\section{Actions on low-dimensional spheres}

In this section we prove that any continuous action of $SL(n,\mathbb{Z})$,
with $n>2$, on the circle and any continuous action of $SL(n,\mathbb{Z})$,
with $n>3$, on $S^{2}$ must be trivial. Since $SL(n,\mathbb{Z})$ is perfect, actions on $S^0$ are trivial.

It follows from the proof of Theorem 1.1 that the actions of $SL(n,\mathbb{Z})$, with $n>2$, on
the circle and the actions of $SL(n,\mathbb{Z})$, with $n>3$, on $S^{2}$
must factor through actions of finite groups, and it is well known that
finite group actions on $S^{1}$ and $S^{2}$ are conjugate to linear actions
of finite subgroups in $O(2)$ and $O(3)$ respectively (see \cite{Edmonds}).
Therefore, it suffices to prove the following theorem.

\begin{theorem}
Every homomorphism from $SL(3+m,\mathbb{Z})$ to $O(2)$ and every
homomorphism from $SL(4+m,\mathbb{Z})$ to $O(3)$, with $m\geq 0$, is trivial.
\end{theorem}

\begin{proof}
First note that Theorem 1.1 implies that the images of the homomorphisms
must necessarily be finite. By a famous result of Milnor, Bass, and Serre,
every homomorphism of $SL(n,\mathbb{Z})$ with finite image factors through $%
SL(n,\mathbb{Z}_{k})$ for some integer $k>1$ (see \cite{Bass}). It is well
known that the group $SL(n,\mathbb{Z}_{k})$ has no abelian quotients. It is
also known, perhaps not well known, that the smallest non-trivial quotient
group of $SL(n,\mathbb{Z}_{k})$ has order at least $2^{n(n-1)/2}\Pi
_{i=2}^{n-1}(2^{i}-1)$ (see \cite{Masur}).

Now consider a homomorphism $h:SL(3+m,\mathbb{Z})\rightarrow O(2)$. Compose $h$
with the determinant map from $O(2)$ to $\mathbb{Z}_{2}$. Every homomorphism
from $SL(3+m,\mathbb{Z})$ to $\mathbb{Z}_{2}$ is trivial because the
congruence subgroups $SL(3+m,\mathbb{Z}_{k})$ have no abelian quotients. This
implies that $h$ factors through $SO(2)$. So $h$ must be trivial because all finite
subgroups of $SO(2)$ are abelian, and again, congruence subgroups $SL(3+m,
\mathbb{Z}_{k})$ have no abelian quotients.

Now consider a homomorphism $g:SL(4+m,\mathbb{Z})\rightarrow O(3)$. The
arguments above show that it suffices to prove that homomorphisms from $SL(4+m,
\mathbb{Z})$ to $SO(3)$ are trivial. All finite subgroups of $SO(3)$ have
been classified---they are all cyclic groups, all dihedral groups, and
all subgroups of the polyhedral groups that have order 12, 24, or 60. The smallest
non-trivial quotient group of $SL(4+m,\mathbb{Z}_{k})$ has order at least 1344
and so subgroups of the polyhedral groups cannot be the images of $g$. Since the congruence subgroups $SL(4+m,\mathbb{Z}_{k})$ have no abelian quotients, it follows that they have no solvable quotients, and therefore, the cyclic and dihedral groups cannot be the images of $g$ also. So every homomorphism from $SL(4+m,\mathbb{Z})$ to $O(3)$ must be trivial.
\end{proof}

In general, continuous actions on higher-dimensional spheres are not
conjugate to linear actions. In \cite{Bing} Bing constructed a continuous $\mathbb{Z}_{2}$ action on $S^{3}$ with Alexander's horned sphere as a fixed-point set. This is impossible for a linear $\mathbb{Z}_{2}$ action and so Bing's action is not conjugate to any linear action. It is conjectured that every smooth action of a compact Lie group on $S^{3}$ is conjugate to a
linear action. However, non-linear smooth actions exist on $S^{n}$ for every 
$n>3$ (see \cite{Gordon}).

It is known that every homomorphism from $SL(n,\mathbb{Z})$ to $GL(r+1,\mathbb{R})$, with $r<n-1$, is trivial. This fact is not well known and the author could not find a reference. A nice proof of this result was shown to the author by Dave Morris, formerly Dave Witte. In particular, this implies that every linear action of $SL(n,\mathbb{Z})$ on $S^{r}$, with $r<n-1$, must be trivial. We conjecture that every continuous action of $SL(n,\mathbb{Z})$ on a $r$-dimensional mod 2 homology sphere, with $r<n-1$, must be trivial.

\section*{Acknowledgments}

The author would like to thank Benson Farb, John Franks, Dave Morris, and
Ronald Dotzel for several useful discussions. The author would also like to thank the referee for carefully reading earlier versions of this paper and for providing many helpful suggestions and comments.

\end{document}